\documentclass[12pt]{article}

\usepackage{graphics}
\usepackage{graphicx}
\usepackage{amssymb}
\usepackage{epsf}

\newcommand{\halmos}{\rule{1ex}{1.4ex}}
\newcommand{\qed}{\hfill \halmos} 
\newcommand{\text}[1]{\hbox{\rm \ #1\ \/}}
\newcommand{\be}[1]{\begin{equation}\label{#1}}
\newcommand{\ee}{\end{equation}}
\newcommand{\bi}{\begin{itemize}}
\newcommand{\ei}{\end{itemize}}
\newcommand{\ben}{\begin{enumerate}}
\newcommand{\een}{\end{enumerate}}

\newcommand{\R}{{\mathbb R}}  
\newcommand{\K}{{\mathcal K}}

\newcommand{\bl}[1]{\begin{Lemma}\label{#1}}
\newcommand{\el}{\qed\end{Lemma}}
\newcommand{\bt}[1]{\begin{Theorem}\label{#1}}
\newcommand{\et}{\end{Theorem}}
\newcommand{\epr}{\end{proof}}
\newcommand{\bpr}{\begin{proof}}
\newenvironment{proof}{\noindent {\em Proof}.\ }{\hspace*{\fill}$\halmos$\medskip}

\newcommand{\beqn}{\begin{eqnarray*}}
\newcommand{\eeqn}{\end{eqnarray*}}

\newtheorem{Theorem}{Theorem}

\newtheorem{Lemma}{Lemma}
\newtheorem{Corollary}{Corollary}

\newcommand{\partes}[1]{\ensuremath{\mathcal{P}(#1)}}
\newcommand{\Norma}[1]{\ensuremath{\parallel \! \! #1 \! \! \parallel}}
\newcommand{\norma}[1]{\ensuremath{\left| #1 \right|}}

\newcommand{\B}{\mathbb{B}}
\newcommand{\Naturals}{\mathbb{N}}
\newcommand{\N}{\mathcal{N}}

\newcommand{\head}[2]{}{}

  \newcommand{\paper}[1]{}

  \newcommand{\stdformat}[1]{}

\title{ \bf Prevalent Behavior of Strongly Order Preserving Semiflows}
\author{Germ\'{a}n A. Enciso \thanks{to whom correspondence should be addressed.  Mathematical Biosciences Institute, Ohio State University, 231 W 18th Ave, Columbus, OH, 43210,
tel. 614 292 6159. Email: genciso@mbi.osu.edu.  This material is
based upon work supported by the National Science Foundation under
Agreement No. 0112050.} \and Morris W. Hirsch\thanks{ Department of
Mathematics, University of California, Berkeley, CA, 94720,
hirsch@math.berkeley.edu} \and Hal L. Smith\thanks{ Department of
Mathematics, Arizona State University, Tempe, AZ 85287,
halsmith@asu.edu, Supported in part by NSF grant DMS 0414270} }

\begin{document}
\maketitle


\begin{abstract}
Classical results in the theory of monotone semiflows give sufficient conditions for the generic solution to converge toward an equilibrium or towards the set of equilibria (quasiconvergence).     In this paper, we provide new formulations of these results in terms of the measure-theoretic notion of prevalence, developed in \cite{Christensen:1972,Yorke:1992}.   For monotone reaction-diffusion systems with Neumann boundary conditions on
convex domains, we show that the set of continuous initial data
corresponding to solutions that converge to a spatially homogeneous equilibrium is
prevalent.
 We also extend a previous generic convergence result to allow its use on Sobolev spaces.  Careful attention is given to the measurability of the various sets involved.
\end{abstract}

\paragraph{Keywords:} strong monotonicity, prevalence, quasiconvergence, reaction-diffusion,
measurability.

\section{Introduction}

The signature results in the theory of monotone dynamics are that
certain dynamic behaviors are generic, for example, convergence to
equilibrium is generic under suitable conditions. In order to be
more precise, some notation is useful but technical definitions will
be deferred to the next section.
Let $\B$ be an ordered
separable Banach space, $X\subset \B$, and consider a semiflow
$\Phi:\R_+\times X\to X$ which is strongly monotone with respect to
a cone $K$ with nonempty interior.
Denote the sets
\[
\begin{array}{l}
B=\{ x\in X\, |\, \mbox{the orbit $O(x)$ has compact closure in $X$} \}  \\
Q=\{ x\in X\, |\, \omega(x)\subseteq E \}   \\
C=\{ x\in X\, |\, \omega(x)=\{e\} \mbox{ for some } e\in E\}.
\end{array}
\]
The elements of $C$ are said to be \emph{convergent} or to have
\emph{convergent solution}, and those of $Q$ are said to be
\emph{quasiconvergent}.

It was established in \cite{Hirsch:1985} that the generic element of
$B$ is quasiconvergent, where  `generic' is made specific in two
different senses:  the topological sense ($B\setminus Q$ is meager
or $Q$ is residual), and the measure theoretic sense ($\mu(B-Q)=0$
for any gaussian measure $\mu$). Later, Smith and Thieme
\cite{Smith:Thieme:Convergence}, motivated by work of Pol\'a\v cik
\cite{Polacik}, provided sufficient conditions for $C$ to contain an
open and dense set.


A drawback of topological genericity is that closed, nowhere dense
subsets of $X$ may still be quite large in terms of measure. In
fact, it is well known that there exists a Cantor subset of $[0,1]$
with positive measure whose complement is open and dense in $[0,1]$.
On the other hand, asking for a set to have measure zero in an
infinite dimensional space $\B$ is difficult to formalize, since
there doesn't exist a measure with the basic properties of the
Lebesgue measure in finite dimensions.  A definition of `sparseness'
that turns out to be very useful in infinite dimensions is that of
prevalence \cite{Christensen:1972,Yorke:1992}:  a set $W\subseteq
\B$ is \emph{shy} if there exists a nontrivial compactly supported
Borel measure $\mu$ on $\B$, such that $\mu(W+x)=0$ for every $x\in
\B$. A set is said to be \emph{prevalent} if its complement is shy.
Given $A\subseteq \B$, we say here that a set $W$ is \emph{prevalent
in $A$} if $A-W$ is shy.  Useful properties of the idea of
prevalence are given in \cite{Yorke:1992}.  Most importantly in the
current paper, a shy set has empty interior, and in finite
dimensions $W$ is shy if and only if $W$ has Lebesgue measure zero.

In this paper, we obtain the counterparts of the genericity results
of Hirsch and Smith and Thieme with prevalence as the notion of
genericity. We show that $Q$ is prevalent in $B$, and that under
additional hypotheses, so is $C$. Genericity in this
measure-theoretic sense seems natural to the theory of monotone
systems. The canonical family of compactly supported Borel measures
for our purposes is given by the uniform measures $\mu_v,\ v>0$,
supported on the segment $J_v=\{tv:0\le t\le 1\}$ joining zero to
the positive vector.   Prevalence with respect to $\mu_v$ is a natural notion of genericity in this context;  for instance, recall Hirsch's result \cite{Hirsch:1988} according to which  $J- Q $ is countable, for any totally ordered arc $J$.


An earlier result \cite{HirschSmith, HirschSmith2} giving sufficient
conditions that $Q$ (and also $C$) contains an open and dense set is
significantly improved. Roughly, this earlier result replaces
certain compactness assumptions on the semiflow by order properties
of the state-space, namely, omega limit sets should have infima or
suprema. We improve it by requiring that some possibly larger space
contains infima or suprema of the limit sets. Our extension
facilitates the application of the theory to partial differential
equations on state spaces that continuously imbed in a space of
continuous functions.
An example of such an application is provided.

If $\Phi$ is $C^1$ and $e\in E$ an equilibrium, define $\rho(e,t)$
to be the spectral radius of the Frechet derivative $D_x\Phi(t,e)$.
We say that $e$ is \emph{linearly stable} if $\rho(e,t)<1$ for all
$t>0$, \emph{linearly unstable} if $\rho(e,t)>1,\ t>0$, and
\emph{neutrally stable} if $\rho(e,t)=1,\ t>0$. Finally, define
$E_s\subseteq E$ to be the set of equilibria that are either
linearly stable or neutrally stable. We show in Theorem~\ref{gen teo main} that the set of initial
data corresponding to trajectories that converge to a point in $E_s$
is prevalent in $X$.
See also Theorem~3.3 and Corollary~3.4 of Pol\'a\v cik
\cite{Polacik02}, where a similar result is discussed whose proof
uses entirely different arguments. Using Theorem~\ref{gen teo main}
and a well known result of Kishimoto and Weinberg
\cite{Kishimoto:1985}, we conclude that for a strongly cooperative
reaction diffusion system of $n$ equations with Neumann boundary
conditions on a convex domain $\Omega$, the set of initial data in
$C(\overline{\Omega},\R^n)$ corresponding to orbits that converge
towards a constant equilibrium is prevalent.

\section{Definitions and Basic Results}

Let $X$ be an ordered metric space with metric $d$ and {\it partial
order} relation $\le$.  We write $x<y$ if $x\le y$ and $x\ne y$.
Given two subsets $A$ and $B$ of $X$, we write $A\le B\ (A<B)$ when
$x\le y\ (x<y)$ holds for each choice of $x\in A$ and $y\in B$. We
assume that the order relation and the topology on $X$ are
compatible in the sense that $x\le y$ whenever $x_n\to x$ and
$y_n\to y$ as $n\to\infty$ and $x_n\le y_n$ for all $n$.  For
$A\subset X$ we write $\bar A$ for the closure of $A$ and
$\hbox{Int} A$ for the interior of $A$. A subset of an ordered space
is {\em unordered} if it does not contain points $x,y$ such that
$x<y$. $A$ is {\em order-convex} if $x\le z\le y$ and $x,y\in A$
implies $z\in A$.

Let $A\subset X$ and let $L=\{x\in X:x\le A\}$ be the (possibly
empty) set of lower bounds for $A$ in $X$. In the usual way, we
define $\inf A:=u$ if $u\in L$ and $L\le u$; $u$ is unique if it
exists. Similarly, $\sup A$ is defined.

The notation $x\ll y$ means that there are open neighborhoods $U,V$
of $x, y$ respectively such that $U\le V$. Equivalently, $(x,y)$
belongs to the interior of the order relation.  The relation $\ll$
is sometimes referred to as the {\em strong ordering}.  We write
$x\ge y$ to mean $y\le x$, and similarly for $>$ and $\gg$.

In most applications, $X$ is a subset of an ordered Banach space
$\B$ having an order cone $\K$. In this case, $x\le y$ if and only
if $y-x\in \K$. If $\K$ has nonempty interior, then $x\ll y$ if and
only if $y-x\in \hbox{Int}\K$. A subset $A$ of $\B$ is {\em
p-convex} if $x<y$ and $x,y\in A$ implies $z=tx+(1-t)y\in A$ for
$0<t<1$.

A {\it semiflow} on $X$ is a continuous map $\Phi:\R_+ \times X\to
X, \:(t,x)\mapsto \Phi_t (x)$ such that:
$$
  \Phi_0(x) = x,  \quad  (\Phi_t \circ \Phi_s)(x)  =\Phi_{t+s}(x)
\qquad(t,s\ge 0,\ x\in X)
$$
The {\it orbit} of $x$ is the set $O(x)=\{\Phi_t(x) : t\ge 0\}$. An
{\it equilibrium} is a point $x$ for which $O(x)=\{x\}$.  The set of
equilibria is denoted by $E$.

The {\it omega limit set} $\omega (x)$ of $x\in X$, defined in the
usual way, is closed and positively invariant. When $\overline
{O(x)}$ is compact, $\omega (x)$ is also nonempty, compact,
invariant, connected, and it attracts $x$.


Let $\Phi$ denote a semiflow in an ordered space $X$. We call $\Phi$
{\it monotone} provided
$$
  \Phi_t(x) \le \Phi_t(y)\ \hbox{whenever}\ x\le y\ \hbox{and}\ t\ge 0.
$$
$\Phi$ is {\it strongly monotone}  if $x<y$ implies that
$\Phi_t(x)\ll \Phi_t(y)$ for all $t>0$ and {\it eventually strongly
monotone} if it is monotone and there exists $t_0\ge 0$ such that
$x<y$ implies that $\Phi_t(x)\ll \Phi_t(y)$ for $t\ge t_0$. $\Phi$
is {\it strongly order-preserving}, SOP for short, if it is monotone
and whenever $x<y$ there exist open subsets $U,V$ of $X$ with $x\in
U$ and $y\in V$ and $t_0\ge 0$ such that
$$
  \Phi_{t_0}(U)\le \Phi_{t_0}(V).$$
Monotonicity of $\Phi$ then implies that $\Phi_t(U)\le \Phi_t(V)$
for all $t\ge t_0$.  Strong monotonicity implies eventual strong
monotonicity which implies SOP. See e.g.
\cite{Smith:Thieme,Smith:monotone,HirschSmith}. The fundamental
properties of an SOP semiflow are stated below.

\bt [Convergence Criterion: \label{CC}  If $\Phi$ is SOP and
$\Phi_T(x)>x$ ($\Phi_T(x)<x$) for some $T>0$ then $\Phi_t(x) \to
p\in E$ as $t\to\infty$. \et

\bt [Nonordering of Limit Sets: \label{NO} Let $\Phi$ be SOP and
$\omega$ be an omega limit set. Then no two points of $\omega$ are
related by $<$. \et

\bt  [Limit Set Dichotomy:\label{LSD} Let $\Phi$ be SOP. If $x<y$
then either
\begin{itemize}
\item [(a)] $\omega(x) < \omega(y)$, or
\item [(b)] $\omega(x) =\omega(y) \subset E$.
\end{itemize}
\et

Smith and Thieme \cite{Smith:Thieme:Convergence} improve part (b) of
the Limit Set Dichotomy to read $\omega(x)=\omega(y)=\{e\}$ for some
$e\in E$ under additional smoothness and strong monotonicity
conditions. For example, this Improved Limit Set Dichotomy holds if
$X\subset Y$ is order convex in the ordered Banach space $Y$ with
cone $Y_+$ having non-empty interior, $\Phi_t(x)$ is $C^1$ in $x$
and its derivative is a compact, strongly positive operator. See
e.g. \cite{Smith:Thieme:Convergence,Smith:monotone,HirschSmith}.

The notion of a shy set, defined in the introduction, was introduced by Yorke et al \cite{Yorke:1992} and
Christensen \cite{Christensen:1972}.
It is a natural generalization to infinite
dimensional spaces of a (Lebesgue) measure zero subset of $\R^n$ in
the sense that a subset of $\R^n$ has measure zero if and only if it
is shy. A countable union of shy sets is shy. Moreover, built into
the definition is translation invariance: if $A$ is shy then so is
$A+x$. Prevalent sets are dense. These and other properties can be
found in \cite{Yorke:1992}.

\section{$C$ is Prevalent in $B$}   \label{gen section C is prevalent}
\head{Prevalence of Convergence}{$C$ is Prevalent in $B$}

Consider an SOP semiflow $\Phi$  defined on a subset $X$ of the
separable Banach space $B$, ordered with respect to a cone $\K$.
Given $v\in B, v\not=0$, define the Borel measure $\mu_v$ on $B$ to
be the uniform measure supported in the set $S_v:=\{ tv \, | 0\leq
t\leq1\}$. That is, $\mu_v(A)=m\{t\in [0,1] \, | tv\in A\}$, where
$m$ is the Lebesgue measure in $[0,1]$.

\begin{Lemma} \label{lemma countable}
Let $W\subseteq X$  be such that $L\cap W$ is countable, for every
straight line $L$ parallel to a positive vector $v> 0$.  Then $W$ is
shy.
\end{Lemma}

\bpr Consider $v>0$ and the uniform measure $\mu_v$.  Let $L=\R v -
x$ for an arbitrary $x\in \B$. Then
\[
(W+x)\cap S_v\subseteq (W+x)\cap \R v=(W\cap L)+x.
\]
Therefore clearly $\mu_v(W+x)=0$, and $W$ is shy with respect to
$\mu_v$. \epr

The proof of Hirsch's generic convergence theorem as stated in terms
of prevalence becomes clear at this point.  See Theorem~4.4 of
\cite{Hirsch:1985}, and \cite{Hirsch:1988}.

\begin{Theorem} \label{gen teo Hirsch}
Let $\B$ be a separable Banach space, and consider a strongly
monotone system defined on $X\subseteq \B$.  Then $Q$ is prevalent
in $B$.

If the Improved Limit Set Dichotomy holds for $\Phi$, then $C$ is
prevalent in $B$.
\end{Theorem}

\bpr  Let $N=B-Q$ be the set of states $x\in B$ such that
$\omega(x)\not\subseteq E$ (see Section~\ref{gen section
measurability} for a proof that this set is measurable). Let
$L\subseteq X$ be a straight line parallel to a vector $v>0$. Note
that if $x,y\in L\cap N$, $x<y$,  then $\omega(x)< \omega(y)$, since
otherwise $\omega(x)=\omega(y)\subseteq E$ by the Limit Set
Dichotomy.


We can apply an argument as in Theorem~7.3 c) of Hirsch
\cite{Hirsch:1988} to conclude that $N\cap L$ is countable: consider
the set $Y=\cup_{x\in N\cap L} \omega(x)$ with the topology
inherited by $\B$.  Since no point in  $\omega(x)$ can bound
$\omega(x)$ from below or above (Theorem~\ref{NO}), no point in
$\omega(x)$ can be the limit of elements in $\omega(y),\ y\not=x$.
Therefore $\omega(x)$ is open in $Y$, for every $x\in N\cap L$.  The
countability of $N\cap L$ follows by the separability of $Y$.

Since $N\cap L$ is countable for every strongly ordered line $L$,
$N$ must be shy by Lemma~\ref{lemma countable}.

If the Improved Limit Set Dichotomy holds for $\Phi$, then we can
argue exactly as above to show that $N=B\setminus C$ is shy. \epr

Define for any set $A\subset X$ the \emph{ strict basin of
attraction}
\[
SB(A):=\{x\in X \, |\, \omega(x)=A\ \mbox{and $x\in B$}\, \}.
\]
Note the difference with the usual basin of attraction of $A$,
$\mathcal{B}(A)=\{x\in X \, |\, \omega(x)\subseteq A\}$.

\begin{Theorem} \label{gen teo main 1}
Let $\B$ be a separable Banach space, and let $X$ be p-convex in
$B$. If $C$ is dense in $B$, then $C$ is prevalent in $B$.
\end{Theorem}

\bpr Let $K=B-C$ be the set of the states $x\in B$ such that
$\omega(x)$ is not a singleton (see Section~\ref{gen section
measurability} for a proof that this set is measurable). We will
show that $K$ is shy with respect to the measure $\mu_v$, for every
$v>0$.  From the assumption that $C$ is dense in $B$, it holds that
$SB(\omega(x))$ has empty interior for every $x\in K$.

Let $L$ be a  straight line in $X$ parallel to a positive vector
$v>0$, $\partes{X}$ denote the set of all subsets of $X$, and
consider the function $\gamma: L\cap K\to
\partes{X}$ defined by $\gamma(x)=\omega(x)$.   Then this function
is injective. Indeed, if $x,y\in L\cap K$, $x<y$, were such that
$\gamma(x)=\gamma(y)=W$, then the strong order preserving property
implies that for any point $u=sx+(1-s)y, \ 0<s<1$ there is a
neighborhood $U$ of $u$ and $t_0\ge 0$ such that $\Phi_t(x)\le
\Phi_t(U)\le \Phi_t(y)$ for $t\ge t_0$ and therefore $\omega(v)=W$
for every $v\in U$ by the Limit Set Dichotomy. As $U$ is a nonempty
open subset of $X$, this implies that $SB(W)$ has nonempty interior
in $X$, a contradiction to the fact that $x$ belongs to $K$.

Note also that, by the Limit Set Dichotomy, the image of $\gamma$ is
an ordered collection of sets: $x<y$ then $\gamma(x)<\gamma(y)$.
Following the same argument as in the proof of Theorem~\ref{gen teo
Hirsch},  it follows that $L\cap K$ is countable.   By
Lemma~\ref{lemma countable}, $K$ is shy with respect to $\mu_v$, for
$v>0$. \epr

\section{$C_s$ is Prevalent in $C$ for Smooth $\Phi$}  \label{gen section stable equilibria}
\head{Prevalence of Convergence}{$C_s$ is Prevalent in $C$}

In this section we assume that $\B$ is a separable Banach space
ordered by a cone $\K$ with nonempty interior and $\Phi$ is a
strongly order preserving semiflow on the p-convex subset $X\subset
\B$. We assume also that for every $t\geq t_0$ the time evolution
operators $\Phi_t$ are compact, (Frechet) $C^1$ and have compact
derivatives (for some fixed $t_0\geq 0$).

We say that an equilibrium point $e\in E$ is \emph{irreducible} if
for some $t=t_e>0$, $\Phi_t'(e)$ is a strongly positive operator
(i.e.\ $x>0$ implies $\Phi_t'(e)x\gg0$). Observe that if
$\Phi_t'(e)$ is strongly positive, so is $\Phi_s'(e), \ s\ge t$. The
point $e$ is said to be \emph{non-irreducible} otherwise. Denote by
$\rho(A)$ the spectral radius of a bounded linear operator $A$ on
$\B$. By the well-known Krein-Rutman theorem \cite{Schaefer}, if $A$
is compact and strongly positive then its spectral radius is a
simple eigenvalue with eigenspace spanned by a positive vector $v\gg
0$; moreover $v$ is the unique eigenvector belonging to $\K$, up to
scalar multiple.

Let
$$
E_s=\{e\in E: \rho(\Phi'_{t_0}(e))\le 1\}
$$
denote the set of ``neutrally stable" equilibria and
$$
C_s=\{x\in X:\omega(x)=\{e\},\ e\in E_s\}
$$
the set of points convergent to a neutrally stable equilibrium. An
equilibrium not in $E_s$ will be called linearly unstable; this
implies that it is unstable in the sense of Lyapunov.

The aim of this section is to provide sufficient conditions for
$C_s$ to be prevalent in $C$ and in $X$.


The following result is well-known.


\begin{Lemma}\label{gen lemma Frechet and equilibria}
Let $T:X\to X$ be a continuous (nonlinear) operator.  Let $e\in X$
be a fixed point of $T$, and assume that the Frechet derivative
$T'(e):T\to T$ exists and is compact.  Assume also that there exists
a sequence $e_1, \ e_2,\ \ldots$ of fixed points of $T$,
$e_k\not=e$, such that $e_k\to e$ as $n\to \infty$.

Then the unit vectors $v_k:=(e_k-e)/\norma{e_k-e}$ have a
subsequence that converges towards a unit vector $w\in \B$, and
$T'(0)w=w$.
\end{Lemma}






\begin{Lemma}
If $(e_k)_{k\in \Naturals}$ is a sequence of equilibria of the
semiflow $\Phi$ such that $e_k< e_{k+1}$ ($e_k< e_{k+1}$) for all
$k$, and if the sequence $(e_k)$ converges towards a irreducible
equilibrium $e\in E$, then $e\in E_s$.
\end{Lemma}

\bpr Let $\tau\geq t_0$ be such that  $\Phi'_\tau(e)$ is strongly
positive, and let $T:=\Phi_\tau$.    Then $T$ satisfies the
hypotheses of Lemma~\ref{gen lemma Frechet and equilibria}, so that
defining $v_k=(e_k-e)/\norma{e_k-e}$, there exists a subsequence
$v_{k_i}$ which converges to a unit vector $w\in \B$.  Furthermore,
$T'(e)w=w$.  From the fact that $e<e_k$ for every $k$, we conclude
that $v_k>0$ and consequently that the unit vector $w>0$ so $w\gg
0$.

By the Krein Rutman theorem, the fact that $T'(e)$ has a positive
eigenvector with eigenvalue 1 implies that in fact $\rho(T'(e))=1$.
Therefore $e\in E_s$, and this concludes the proof.

The case $e_{k+1}< e_k$ for every $k\in \Naturals$ can be treated
similarly. \epr

We introduce the property (P):

\begin{description}
\item[(P)] Every set of equilibria $\hat{E}\subseteq E$ which is totally
ordered by $<$ has at most enumerably many non-irreducible points.
\end{description}

For instance, this condition holds if all equilibria in $X$ are
irreducible (see condition (S) in \cite{Smith:monotone}, p.\ 19). It
also holds if every totally ordered subset of $E$ has at most
enumerably many points.


\begin{Lemma} \label{gen lemma discrete}
Let property (P) be satisfied.   If $\hat{E}\subseteq E$ is totally
ordered by $<$, and if every element of $\hat{E}$ is linearly
unstable, then $\hat{E}$ is countable.
\end{Lemma}

\bpr Suppose that $\hat{E}$ is not countable.   Then the set
$\tilde{E}\subseteq \hat{E}$ of irreducible elements in $\hat{E}$ is
also uncountable, by property (P).      Let $e\in \tilde{E}$ be an
accumulation point of $\tilde{E}$, which exists by the separability
of the Banach space $\B$ (otherwise $\B$ would contain an
uncountable set of pairwise disjoint open balls).     Then there
exists a monotone sequence of elements in $\tilde{E}$ which
converges towards $e$.    By the previous lemma it holds that $e\in
E_s$, contradicting $e\in \tilde{E}$. \epr

\begin{Lemma} \label{gen lemma unstable}
If $e\in E\setminus E_s$, then $SB(e)$ is unordered and hence shy.
\end{Lemma}

\bpr The same argument as Lemma~2.1 in
\cite{Smith:Thieme:coexistence} shows that $SB(e)$ is unordered.
This implies that it is shy with respect to $\mu_v$ for any $v>0$,
by Lemma~\ref{lemma countable}, since any line parallel to $v$ meets
$SB(e)$ at most once.
 \epr

Our next result is similar to Theorem~4.4 in \cite{Hirsch:1985} in
finite dimensions, and to a lesser extent to Theorem~10.1 in
\cite{Hirsch:1988} but it drops the assumptions of finiteness or
discreteness for the set $E$.

\begin{Theorem} \label{gen teo main 2}
In addition to the assumptions of this section, let property (P) be
satisfied.  Then  $C_s$ is prevalent in $C$.
\end{Theorem}

\bpr We follow a very similar argument as in the proof of
Theorem~\ref{gen teo Hirsch}.    Let $N=C-C_s$ be the set of $x\in
C$ such that $\omega(x)$ is a linearly unstable equilibrium. It will
be shown in Section~\ref{gen section measurability} that this set is
Borel.  Let $v>0$ and let $L$ be a line parallel to $v$. Then we can
define the function $\sigma:L\cap N\to X$ by $\sigma(x)=\lim_{t\to
\infty} \Phi(t,x)$.   If $x_1,x_2\in L\cap N$, $x_1< x_2$, then
necessarily $\sigma(x_1)<\sigma(x_2)$ by the Limit Set Dichotomy
since $SB(\omega(x_1))$ is unordered by Lemma~\ref{gen lemma
unstable}. Thus $\sigma$ is injective. As $\hat{E}=\mbox{range }
\sigma$ is totally ordered, it is countable by Lemma~\ref{gen lemma
discrete}, and so is $L\cap N$ by injectivity. By Lemma~\ref{lemma
countable}, $N$ is shy. \epr

See also Theorem~4.4 and Theorem~4.1 of \cite{Hirsch:1985}.

\begin{Theorem} \label{gen teo main}
In addition to the assumptions of this section, suppose $B=X$, $X$
is order convex, every equilibrium is irreducible, and that $\Phi$
is eventually strongly monotone. Then $C_s$ is prevalent in $X$.
\end{Theorem}

\bpr In this case, the Improved Limit Set Dichotomy holds so $C$ is
prevalent in $X$ by Theorem~\ref{gen teo Hirsch}. As $C_s$ is
prevalent in $C$ by Theorem~\ref{gen teo main 2}, the result follows
since $X\setminus C_s=(X\setminus C)\cup (C\setminus C_s)$ is the
union of two subsets, each shy relative to the same $\mu_v, \ v>0$.
\epr

\section{$Q$ Contains an Open and Dense set}

Let $\Phi$ be an SOP semiflow on the ordered metric space $X$,
having  compact orbit closures. In this section we improve a result
in \cite{HirschSmith2} by weakening the conditions for $Q$ to
contain an open and dense set. We introduce the following
hypothesis:

\medskip
\noindent {\bf (K)} $X\subset Z$ where $Z$ is an ordered metric
space with order relation $\le_Z$, the inclusion $i:X\to Z$ is
continuous and $X$ inherits its order relation from $Z$.  $\Phi$
extends to a mapping $\Psi: \R^+ \times Z\to Z$, not necessarily
continuous, where
\begin{enumerate}

\item[(a)]$\Psi|_{\R^+ \times X}=\Phi$, \\

\item[(b)] $\Psi$ is monotone on $Z$.\\

\item[(c)] For every $z\in Z$, there exists $t_z\ge 0$ such that
$\Psi_{t_z}(z)\in X$.

\end{enumerate}

Observe that if $z\in Z$, then $\Psi_t(z)\in X$ and
$\Psi_t(z)=\Phi_t(z)$ for $t\ge t_z$. Consequently, the omega limit
set of the orbit through $z$ exists in both $X$ and in $Z$ and they
agree.

A point $x\in X$ is {\em doubly accessible from below}
(respectively, above) if in every neighborhood of $x$ there exist
$f, g$ with $f<g<x$ (respectively, $x<f<g$).

For $p\in E$ define $C(p):=\{z\in X: \omega(z)=\{p\}\}$.  Note that
$C=\bigcup_{p\in E}C(p)$. All topological properties used hereafter
are relative to the space $X$.

\begin{Lemma}  \label{NGQlem}
Let (K) hold. Suppose $x\in X\setminus Q$ and $a=\inf \omega(x)\in
Z$ exists. Then $\omega(a)=\{p\}$ where $p\in X$ satisfies
$p<\omega(x)$, and $x\in \overline{\hbox{\em Int\,}C(p)}$ provided
$x$ is doubly accessible from below.
\end{Lemma}
\bpr Fix an arbitrary neighborhood $M$ of $x$. Note that $a<_Z
\omega (x)$ because $\omega (x)\subset X$ is unordered
(Theorem~\ref{NO}). By invariance of $\omega (x)$ and (K) we have
$\Psi_t a\le_Z \Psi_t \omega(x)=\Phi_t\omega(x)=\omega(x)$, hence
$\Psi_t a \le_Z a, t\ge 0$. It follows from (K) that $\Phi_t(w)\le
w:=\Psi_{t_a}(a)$ for $t\ge 0$ and therefore the Convergence
Criterion Theorem  implies that $\omega (a)$ is an equilibrium $p\in
X$ with $p\le a$. Because $p<\omega(x)$, SOP yields a neighborhood
$N$ of $\omega (x)$ and $s\ge 0$ such that $p\le \Phi_t N$ for all
$t\ge s$. Choose $r\ge 0$ with $\Phi_t x\in N$ for $t\ge r$. Then
$p\le\Phi_{t}x$ if $t\ge r+s$. The set $V:= (\Phi_{r+s})^{-1}
(N)\cap M$ is a neighborhood of $x$ in $M$ with the property that
$p\le \Phi_t V$ for all $t\ge r+2s$. Hence:
\begin{equation}        \label{eq:ngq0}
u\in V \Rightarrow  p \le\omega(u)
\end{equation}

Now assume $x$ doubly accessible from below and fix $y_1,y\in V$
with $y_1<y<x$. By the Limit Set Dichotomy $\omega(y) <\omega(x)$,
because $\omega(x)\not\subset E$. By SOP we fix a neighborhood
$U\subset V$ of $y_1$ and $t_0>0$ such that $\Phi_{t_0} u \le
\Phi_{t_0} y$ for all $u\in U$.  The Limit Set Dichotomy implies
$\omega(u)=\omega(y)$ or $\omega(u)<\omega(y)$; as
$\omega(y)<\omega (x)$, we therefore have:

\begin{equation}        \label{eq:ngq1}
u\in U \Rightarrow  \omega(u)< \omega(x)
\end{equation}
For all $u\in U$,  (\ref{eq:ngq1}) implies $\omega(u)\le
\omega(a)=\{p\}$, while (\ref{eq:ngq0}) entails $p\le \omega(u)$.
Hence $U \subset C(p)\cap M$, and  the conclusion follows. \epr

An analogous result holds if ``$a=\inf\omega(x)\in Z$ exists" is
replaced by ``$b=\sup\omega(x)\in Z$ exists", in which case
$\omega(b)=\{q\}$ where $q>\omega(x)$. Furthermore, the conclusion
$x\in \overline{\hbox{\em Int\,}C(p)}$ holds provided $x$ is doubly
accessible from above.

We introduce an additional condition on the semiflow $\Phi$ similar
to the one in \cite{HirschSmith2}:

\medskip
\noindent {\bf (L) } Either every omega limit set $\omega(x),\ x\in
X$, has an infimum in $Z$ and the set of points that are  doubly
accessible from below has dense interior in $X$, or every omega
limit set has a supremum in $Z$ and the set of points that are
doubly accessible from above has dense interior in $X$.
\medskip

\begin{Theorem}  \label{NGQ}
Let $\Phi$ be an SOP semiflow on the ordered metric space $X$,
having  compact orbit closures, and satisfying axioms {\rm (L)} and
{\rm(K)}. Then $X\setminus Q \subset \overline{\hbox{\em Int\,}C}$,
and $\hbox{\em Int\,}Q$ is dense.
\end{Theorem}

\bpr To fix ideas we assume the first alternative in (L), the other
case being similar.  Let $X_0$ denote a dense open set of points
doubly accessible from below. Lemma \ref{NGQlem} implies $
  X_0\subset Q\cup \overline{\hbox{\rm Int\,}C} \subset Q\cup
 \overline{\hbox{\rm Int\,}Q},
$ hence the open set $X_0\setminus \overline{\hbox{\rm Int\,}Q}$
lies in $Q$.  This prove $X_0\setminus \overline{\hbox{\rm Int\,}Q}
\subset \hbox{\rm Int\,}Q$, so $X_0\setminus \overline{\hbox{\rm
Int\,}Q}=\emptyset$. Therefore $\overline{\hbox{\rm  Int\,}Q}\supset
X_0$, hence $\overline{\hbox{\rm Int\,}Q}\supset\overline{X_0}=X$.
\epr

Axiom (L) is a restriction on both the space $X$ (order and
topology) and the semiflow (limit sets). If $X$ continuously embeds
in $Z=C(A,\R)$, the Banach space of continuous functions on a
compact set $A$ with the usual ordering, $\Phi$ extends to a
monotone mapping on $Z$ with the smoothing property (c), then axiom
(L) holds. This is true because every compact subset of $C(A,\R)$
has a supremum and infimum (see Schaefer \cite{Schaefer}, Chapt. II,
Prop. 7.6). In particular, $X$ may be a Euclidean space $\R^n$, a
H\"{o}lder space $C^{k+\alpha}(\Omega,\R^n),\ 0\le\alpha<1,
k=0,1,2,\cdots,$ for $\Omega$ a compact smooth domain in $R^m$, or a
Sobolev space $H^{k,p}(\Omega)$ for $k-\frac{n}{p}\ge 0$ where the
usual functional ordering is assumed. These cases cover ordinary,
delay, and parabolic partial differential equations under suitable
hypotheses.

Theorem~\ref{NGQ} extends the corresponding result in
\cite{HirschSmith2}, where it was assumed that $Z=X$, by allowing
$X$ to be imbedded in a larger space $Z$ in which it is more likely
that omega limit sets have infima and suprema. This extension is
important for partial differential equations for the reasons
mentioned above.

We show how Theorem~\ref{NGQ} can be used to improve Theorem 6.17 of
\cite{HirschSmith} concerning the system of $m$ reaction diffusion
equations given by
\begin{eqnarray*}
\frac{\partial u_i}{\partial t}&=&A_iu_i+f_i(x,u),\ x\in \Omega,\
t>0\\
B_iu_i&=&0,\ x\in \partial \Omega,\ t>0
\end{eqnarray*}
where $A_i$ are uniformly elliptic second order differential
operators and $B_i$ are boundary operators of Dirichlet, Robin, of
Neumann type and $f$ is cooperative and irreducible in the sense
that  $f(x,u)$ is $C^1$ in $u$ and $\partial f_i/\partial u_j \ge 0$
for all $i\ne j$. In addition, there exists $\bar x\in \Omega$ such
that the $m\times m$ Jacobian matrix  $[\partial f_i/\partial
u_j(t,\bar x,u)]$ is irreducible for all $u$. These conditions could
be formulated with respect to an alternative orthant order with no
change in conclusions.

Let  $\Gamma\subset \R^m$ be a rectangle, i.e., product of $m$
nontrivial intervals, and $k=0,1$ define
$$
X^k_\Gamma :=\{u\in C^k_B (\overline\Omega, \R^m):u
   (\overline\Omega)\subset \Gamma\},\quad X_\Gamma :=\{u\in L^p(\Omega, \R^m):u
   (\Omega)\subset \Gamma\}
$$
See \cite{HirschSmith} for details on the notation; the subscript
$B$ indicates the boundary conditions are accommodated. We assume
that the system above generates a semiflow $\Phi$ on $X_\Gamma$ and
semiflows $\Phi^k,\ k=0,1$ on $X^k_\Gamma$. See \cite{HirschSmith}
for such conditions.

Finally, assuming that these semiflows have compact orbit closures,
it is observed in \cite{HirschSmith} that:

\begin{itemize}

\item $X^1_\Gamma$ is dense in $X^0_\Gamma$ and in $X_\Gamma$.

\item $\Phi$ and $\Phi^0$ agree on $X^0_\Gamma$, and $\Phi,
\Phi^0$ and $\Phi^1$ agree on $X^1_\Gamma$

\item $\Phi_t$ (respectively, $\Phi^0_t$) maps $X_\Gamma$
 (respectively, $X^0_\Gamma$) continuously into
 $X^1_\Gamma$ for $t>0$

\item  $\Phi$, $\Phi^1$ and $\Phi^0$ have the same omega limit
sets, compact attractors and equilibria.

\end{itemize}

Theorem 6.17 in \cite{HirschSmith} concludes, among other things,
that the set of quasiconvergent points for each of the semiflows is
residual in the appropriate space. In particular, $Q(\Phi^k)$ is
residual in $X^k_\Gamma$ for $k=0,1$. In fact, $Q(\Phi^k)$ is open
and dense in $X^k_\Gamma$ by Theorem~\ref{NGQ}. To see this, we need
only note that $X^1_\Gamma$ imbeds continuously in  $Z=X^0_\Gamma$,
that axiom {\rm(K)} holds by virtue of the properties itemized
above, and that axiom {\rm (L)} holds in $Z$ for the reasons noted
following the proof of Theorem~\ref{NGQ}.

\section{Reaction-Diffusion Systems with No-Flux Conditions on a Convex Domain}  \label{gen section rd}
\head{Prevalence of Convergence}{Reaction-Diffusion Systems}

Consider a reaction-diffusion system \begin{eqnarray}\label{gen eq
rd}
u_t&=&D \Delta u + f(u),\ x\in \Omega\nonumber\\
\frac{\partial u}{\partial n}&=&0, \ x\in \partial\Omega\\
u(0,x)&=&u_0(x),\ x\in \overline{\Omega}\nonumber
\end{eqnarray}
We assume that the state space $C(\overline{\Omega},\R^n)$.
Kishimoto and Weinberger \cite{Kishimoto:1985} showed that if
$\Omega$ is a convex domain, and assuming that $\partial f_i /
\partial u_j >0 $ for all $i\not=j$, then any nonconstant
equilibrium $\overline{u}$ is linearly unstable.  A careful reading
of the proof in that paper will show that in fact it is sufficient
that $\partial f_i /
\partial u_j \geq 0 $ for all $i\not=j$ and that the Jacobian matrix is irreducible.
We call the vector field $f$ cooperative and irreducible if this is
the case. By making a linear change of variables, we may extend the
following result to any system which is monotone with respect to one
of the other orthants $\K$. See \cite{Smith:monotone}.

Equilibria of (\ref{gen eq rd}), solutions of the associated
elliptic boundary value problem, are known for having multiple and
sometimes unexpected solutions. Not only is it possible for a
strongly monotone reaction diffusion system to have several
spatially nonhomogeneous equilibria, but it is in fact possible that
there is a continuum of them \cite{Enciso:continuum}.  The following
application of Theorem~\ref{gen teo main} shows that the generic
solution converges to a uniform (constant) solution.

\begin{Theorem}  \label{gen teo rd conclusion}
Consider a $C^1$ finite dimensional system \be{gen equation
diffusion-finite} \dot{x}=f(x) \ee which is cooperative and
irreducible  and assume that all initial value problems have bounded
solutions for $t\ge 0$. If $\Omega$ is convex, then the set of
initial conditions $u_0\in C(\overline{\Omega},\R^n)$ corresponding
to solutions of (\ref{gen eq rd}) that converge towards a uniform
equilibrium is prevalent in $C(\overline{\Omega},\R^n)$.
\end{Theorem}

\bpr We need to show that all the general assumptions of the
previous section are satisfied, as well as the hypotheses of
Theorem~\ref{gen teo main}.  Clearly $X=C(\overline{\Omega},\R^n)$
is a separable Banach space under the uniform norm with cone given
by $C(\overline{\Omega},\R_+^n)$. The fact that the time evolution
operators generate a semiflow of compact operators with compact
derivatives on $X$ is well known in the literature; see for instance
\cite{Lunardi,Pazy,HirschSmith,Polacik02}.   The fact that $B=X$
follows from comparison with solutions of the ordinary differential
equations (\ref{gen equation diffusion-finite}); see e.g.
Theorem~7.3.1 in \cite{Smith:monotone}. To see that the system
(\ref{gen eq rd}) has no non-irreducible equilibria, let $\hat{u}$
be an equilibrium of the system, and recall that the linearization
around $\hat{u}$ is of the form
\[
u_t=D\Delta u + M(x)u,
\]
together with Neumann boundary conditions, where  $M(x)=\partial f /
\partial u (\hat{u}(x))$. According to Theorem~7.4.1 of
\cite{Smith:monotone}, to prove that this system is strongly
monotone it is enough to verify that the associated
finite-dimensional system with no diffusion is monotone for every
fixed value of $x\in \Omega$, and strongly monotone for at least one
value of $x$.  This therefore follows from the irreducibility
assumption on the linearizations of (\ref{gen equation
diffusion-finite}).

By Theorem~\ref{gen teo main}, $C_s$ is prevalent in $B=X$.  But by
the main theorem in \cite{Kishimoto:1985}, any initial condition in
$C_s$ has a solution which converges towards an equilibrium which is
uniform in space. This completes the proof. \epr

\section{Appendix: Measurability}  \label{gen section measurability}
\head{Prevalence of Convergence}{Regarding Measurability}

It is important to observe that in order to apply measure-theoretic
arguments on Theorems~\ref{gen teo Hirsch}, \ref{gen teo main 1},
and \ref{gen teo main 2}, one needs to prove first that the sets
involved in each result are Borel measurable.  This is carried out
in the present section, where we assume throughout that $X$ is a Borel subset of a \emph{separable} Banach space $\B$,
that $\Phi:\hbox{Dom}\ \Phi\to X$ is a continuous local semiflow defined
on the open subset $\hbox{Dom}\ \Phi$ of $R_+\times X$ containing
$\{0\}\times X$. For each $x\in X,\ \{t\ge 0: (t,x)\in \hbox{Dom}\
\Phi \}=[0,\sigma_x)$. The set
$$
Ext=\{x\in X: \sigma_x=+\infty\}=\cap_{q} \{x\in X: \sigma_x>q\},
$$
where the intersection is taken over all positive rational $q$, is
Borel since it is the countable intersection of open subsets of $X$.
Therefore, we may as well rename $X=Ext$ and consider the global
semiflow $\Phi: \R_+\times X\to X$ where $X$ is Borel.  Given $p\in
\B$ and $r>0$, let $B_r(p)=\{x\in \B:|x-p|<r\}$.

Let $D\subseteq X$ be a closed set in $X$ and  $r\in R_+$, and
consider the set
\[
W(D,r)=\{x\in X\, |\, \Phi_t(x)\in D,\ \mbox{for all } t\geq
r\}=\bigcap_{q\in Rat,\ q>r} \Phi_q^{-1}(D),
\]
where $Rat$ denotes the rationals. The equality holds  from the
continuity of $\Phi$. Since each operator $\Phi_q$ is continuous,
$W(D,r)$ is a Borel measurable set.

In the following we assume only that $\Phi$ is a continuous semiflow
on the closed set $X$.

\begin{Lemma} \label{lemma B}  If $X$ is closed in $\B$, then the set $B$ of the elements $x\in X$ with
precompact orbit is Borel measurable.
\end{Lemma}

\bpr Note that $B$ is the set of $x\in X$ such that $O(x)$ is totally
bounded, and that a set $S$ is totally bounded if and only if for every $\epsilon>0$ there exists a finite collection of \emph{closed} balls of radius less than $\epsilon$, whose union contains $S$.   Let $\{p_i\}_{i\in
\mathbb{N}}$ be a countable dense set in $\B$ and let $\mathfrak{F}$
be the family of all finite subsets of $\mathbb{N}$. Then
$\mathfrak{F}$ is countable and
\[
B=\bigcap_{n\in \mathbb{N}}\bigcup_{F\in \mathfrak{F}} W(S_{F,n},0),
\]
where $S_{F,n}=\cup_{i\in F} \overline{B}_{1/n}(p_i)$.  It is easy to see from this expression that $B$ must be Borel measurable.


\epr

\begin{Lemma}
Let $D\subseteq X$ be a closed set, and let $C(D)=\{x\in B\,|\,
\omega(x)\subseteq D\}$.  Then $C(D)$ is Borel measurable.
\end{Lemma}

\bpr
Given a set $A\subseteq X$ and $\epsilon>0$, let
\[
A_\epsilon=\{x\in X\,|\, d(A,x)\leq\epsilon \},
\]
which is a closed set by continuity of the function $d(\cdot, A)$.
Then we can write
\[
\{x\in B\, |\, \lim_{t\to \infty} d(\Phi_t(x),A)=0\}=\bigcap_{m\in
\Naturals}\bigcup_{k\in \Naturals}W(A_{\frac{1}{m}},k).
\]
Finally, note that for any closed set $D\subseteq X$ and for any
$x\in B$, it holds that
\[
\omega(x)\subseteq D \Leftrightarrow \lim_{t\to \infty}
d(\Phi_t(x),D)=0.
\]
The first statement follows.
\epr

\begin{Corollary} The set $Q$ of quasiconvergent elements is Borel measurable.
\end{Corollary}
\bpr The proof follows immediately from the above result by noting that $E$ is a closed set.
\epr

It follows that the set $N=B-Q$ involved in the proof of
Theorem~\ref{gen teo Hirsch} is measurable, by the previous
Corollary and Lemma~\ref{lemma B}.

\begin{Lemma}  The set $C$ of convergent elements is Borel measurable.
\end{Lemma}

\bpr Let $(b_i)$ be a dense enumerable collection of elements of
$X$.  Then the statement follows from the equation
\[
C=\bigcap_{k=1}^\infty \bigcup_{i=1,r=1}^\infty W(\overline{B}_{\frac{1}{k}}(b_i),r).
\]
To see this, let first $x\in C$.  Note that for any fixed $k$, there
exists some $b_i$ within $1/(2k)$ of $\omega(x)$, and that therefore
$x\in W(\overline{B}_{\frac{1}{k}}(b_i),r)$ for some large enough
$r$.  Therefore for every fixed $k$, $C$ is contained in the union
of the RHS, and thus one direction is proven.  Conversely, let $x$
be in the RHS term.  For every $k$, let $a_k,r_k$ be such that $x\in
W(\overline{B}_{\frac{1}{k}}(a_k),r_k)$; such sequences exist by
hypothesis.  Note that for $k_1\not=k_2$, it must hold
\[
W(\overline{B}_{\frac{1}{k_1}}(a_{k_1}),r_{k_1})\cap
W(\overline{B}_{\frac{1}{k_2}}(a_{k_2}),r_{k_2}) \not=\emptyset
\]
In particular the sequence $(a_k)$ is Cauchy, and it therefore
converges towards a point $a\in \B$.  Given $\epsilon>0$, let $k$
large enough that $1/k < \epsilon/2$ and that $a_k$ lies within
$\epsilon/2$ of $a$.  By definition, $\norma{x(t)-a}<\epsilon$ for
$t\geq r_k$;  we conclude that $x(t)\to a$, and therefore that $a\in
X$ and $x\in C$. \epr

\begin{Lemma}
Assume that the time evolution operator is continuously
differentiable. Then the set $C_s$ is Borel measurable.
\end{Lemma}

\bpr Note first that the spectral radius function $\rho(T)$, though
not a continuous function of the linear bounded operator
$T:L(\B,\B)\to \R$ (see Kato \cite{Kato}), is nevertheless a
measurable function.  To see this, simply write it as the pointwise
limit of continuous functions as $\rho(T)=\lim_n \Norma{T^n}^{1/n}$.
Fix now $t>t_0$, and define $\beta:X\to \R$,
$\beta(z):=\rho(\Phi_t'(z))$.  Since $z\to \Phi_t'(z)$ is a
continuous function by hypothesis, it follows that $\beta$ is
measurable.

The next step is to note that the function $z\to \omega(z)$ (defined
on $C$) is also measurable.  To see this, write it as the pointwise
limit of the continuous functions $\omega(z)=\lim_n \Phi_n(z)$.
Thus, the function $z\to \beta(\omega(z))$ is itself measurable. But
\[
C_s=\{z\in C\, |\, \beta(\omega(z))\leq 1\}, \] and the proof is
complete. \epr

Note that the continuous differentiability of the time evolution
operators was only used to show that $\beta$ is measurable; it would
be sufficient to assume $z\to \Phi_t'(z)$ to be measurable, which
should be satisfied in very large generality.

\begin{Lemma}  Let $A\subseteq B$ be compact.  Then $SB(A)$ is Borel measurable.
\end{Lemma}

\bpr For every $\epsilon>0$, there exists a finite collection
$R_\epsilon$ of open balls of radius $\epsilon$, such that i) each
ball intersects $A$, and ii) the union of all balls contains $A$.
Let $R=\cup_{n\in \N} R_{1/n}$.   Then the set $\bigcup_{V\in R}
W(V^C,0)$ consists of the vectors $x\in X$ such that $a\not\in
\omega(x)$ for some $a\in A$.  Consequently,
\[
SB(A) =  C(A) - \bigcup_{V\in R} W(V^C,0),
\]
and this set is also Borel measurable. \epr

\paragraph{Acknowledgements}
The first author would like to thank Liming Wang for her help with
the generalization of the result in \cite{Kishimoto:1985}, and
Eduardo Sontag for several useful discussions.

\end{document}